\documentclass[12pt,a4paper]{amsart}
\usepackage[top=1.15in, bottom=1.15in, left=1.17in, right=1.17in]{geometry}
\usepackage{amssymb}
\usepackage{pdflscape}
\usepackage{amscd}
\usepackage{fancyhdr}

\newtheorem{theorem}{Theorem}[section]
\newtheorem{corollary}[theorem]{Corollary}
\newtheorem{lemma}[theorem]{Lemma}
\newtheorem{proposition}[theorem]{Proposition}

\usepackage[usenames,dvipsnames]{xcolor} 
\usepackage{tikz, tikz-3dplot, pgfplots}
\usepackage{tkz-graph}
\usepackage{tikz-cd}
\usetikzlibrary{calc}

\title{Momentum map images of representation spaces of quivers}
\author{Maria Bertozzi, Markus Reineke}
\begin{document}
\begin{abstract}{We consider the base change action on real or complex representation spaces of quivers and the associated momentum map for a maximal compact subgroup of the base change group, as introduced by A. King. We give an explicit description of the momentum map image in terms of recursively defined inequalities on eigenvalues of Hermitian operators. Moreover, we characterize when the momentum map image is maximal possible, respectively of positive volume.}\end{abstract}
\maketitle
\parindent0pt

\section{Introduction}

Given a Hamiltonian group action of a compact Lie group on a symplectic manifold, it is an interesting problem to determine the image of the associated momentum map to the dual of the Lie algebra of the group. By famous results of Atiyah, Guillemin-Sternberg, Kirwan and Sjamaar \cite{A,GS,Kir,Sj}, this image intersects the positive Weyl chamber in a convex polyhedral cone.\\[1ex]
In this note, we consider the base change action on real or complex representation spaces of quivers, and the momentum map for a maximal compact subgroup, as introduced in \cite{Kin}. We give an explicit description of the momentum map image in terms of recursively defined inequalities on eigenvalues of Hermitian operators, see Theorem \ref{mainintro}. This result is achieved using a combination of standard techniques from the representation theory of quivers due to Schofield \cite{Sch}, and ideas used in work of Derksen-Weyman for proving Horn's conjecture using quivers \cite{DW}, and in work of Hausel-Letellier-Rodriguez Villegas on the Kac conjectures \cite{HLV}.\\[1ex]
Both the result, and the methods for its proof, are similar to the approach of \cite{BVW1,BVW2} determining the highest weights of polynomial functions on representation spaces and the Schubert positions of subrepresentations of quiver representations. In fact, our main result Theorem \ref{mainintro} is not new, and can be obtained by a direct translation of \cite[Theorem 5.1]{BVW2}. However, the authors feel that a short and direct description of the momentum map image in terms of symplectic geometry, as well as a concise derivation of this description, are desirable.\\[1ex]
In Section \ref{momentum}, we introduce the setup of the base change action on complex representation spaces and formulate the main result, without reference to notions of the representation theory of quivers. Representation-theoretic techniques, in particular the notion of general subrepresentations, are recollected in Section \ref{general}, following \cite{Ch,Kin,Sch}. The main technique for determining the momentum map image, namely translation to the so-called leg-extended quiver, is introduced in Section \ref{leg} (relying on a lemma from \cite{CG}), and leads to the representation-theoretic description Corollary \ref{cor51}. We perform several technical reductions in Section \ref{reductions} to make this criterion effective, culminating in the main result Theorem \ref{mainintro}. We illustrate our main result with a few  examples in Section \ref{examples}, and apply results of \cite{LBP,Sch} to derive sufficient criteria for the momentum map image to be maximal, respectively of positive volume. In Section \ref{real} we consider the case of real representation spaces using \cite{OS}.\\[2ex] 
{\bf Acknowledgments:} The authors would like to thank P. Heinzer for posing the question whether an explicit description of the momentum map image can be obtained, and P. Albers and W. Crawley-Boevey for helpful remarks on the main result. The authors are indebted to V. Baldoni, M. Vergne and M. Walter for pointing out errors in an earlier version, and for thorough discussions on the comparison of the present work to \cite{BVW1,BVW2}. This work was supported by the DFG-CRC/TRR 191 Symplectic Structures in Geometry, Algebra and Dynamics.

\section{The momentum map}\label{momentum}

Let $Q$ be a finite quiver, that is, an oriented graph with a finite set of vertices $Q_0$ and finitely many arrows written $\alpha:i\rightarrow j$. Let $\mathbf{d}=(d_i)_i\in\mathbb{N}Q_0$ be a dimension vector for $Q$, and fix complex vector spaces $V_i$ of dimension $d_i$ for all $i\in Q_0$, respectively. Let $$R_\mathbf{d}(Q)=\bigoplus_{\alpha:i\rightarrow j}{\rm Hom}_\mathbb{C}(V_i,V_j)$$
be the space of representations of $Q$ on the vector spaces $V_i$. The group
$${\rm GL}_\mathbf{d}=\prod_{i\in Q_0}{\rm GL}(V_i)$$
acts on $R_\mathbf{d}(Q)$ via the base change action
$$(g_i)_i\cdot(f_\alpha)_\alpha=(g_j\circ f_\alpha\circ g_i^{-1})_{\alpha:i\rightarrow j}.$$
The orbits for this action are, by definition, precisely the isomorphism classes of complex representations of $Q$ of dimension vector $\mathbf{d}$.\\[1ex]
We follow King \cite[Section 6]{Kin} in describing an associated momentum map. A choice of a Hermitian form on every vector space $V_i$ yields a Hermitian form on $R_\mathbf{d}(Q)$ via
$$((f_\alpha),(g_\alpha))=\sum_\alpha{\rm trace}(f_\alpha\circ g_\alpha^*).$$
 
The product of the associated unitary groups $U_\mathbf{d}=\prod_{i\in Q_0}U(V_i)$ is a maximal compact subgroup of ${\rm GL}_\mathbf{d}$ preserving the Hermitian form on $R_\mathbf{d}(Q)$. Denoting by $\mathfrak{u}_\mathbf{d}$ the Lie algebra of $U_\mathbf{d}$, the infinitesimal action of $i\mathfrak{u}_\mathbf{d}$ on $R_\mathbf{d}(Q)$ is given by $$(A_i)_i\cdot(f_\alpha)_\alpha=(A_j\circ f_\alpha-f_\alpha\circ A_i)_{\alpha:i\rightarrow j}.$$

A momentum map $\mu:R_\mathbf{d}(Q)\rightarrow (i\mathfrak{u}_\mathbf{d})^*$ for the action of $U_\mathbf{d}$ on $R_\mathbf{d}(Q)$ is then provided by $$\mu((f_\alpha))((A_i))=((A_i)\cdot(f_\alpha),(f_\alpha)).$$

Identifying $(i\mathfrak{u}_\mathbf{d})^*$ with the space $\bigoplus_{i\in Q_0}{\rm Herm}(V_i)$ of tuples of Hermitian operators, again via the standard trace form, we get the following description of the momentum map $\mu=\mu_Q:R_\mathbf{d}(Q)\rightarrow\bigoplus_{i\in Q_0}{\rm Herm}(V_i)$:
$$\mu((f_\alpha)_\alpha)=(\sum_{\alpha:~\rightarrow i}f_\alpha \circ f_\alpha^*-\sum_{\alpha:i\rightarrow ~}f_\alpha^*\circ f_\alpha)_i.$$

{\bf Remarks:}
\begin{enumerate}
\item The target $(i\mathfrak{u}_\mathbf{d})^*$ of the momentum map is chosen to be consistent with \cite{Kin}. Considering a different dualization to $\bigoplus_{i\in Q_0}{\rm Herm}(V_i)$, one could also choose the more common target $\mathfrak{u}_\mathbf{d}^*$, which is then a momentum map for the Hamiltonian action of $U_\mathbf{d}$ on $R_\mathbf{d}(Q)$, where the latter is endowed with the symplectic structure $\omega=2{\rm Im}(\_,\_)$ derived from the Hermitian structure.
\item Our setup should not be confused with the momentum map for the algebraic symplectic ${\rm GL}_{\mathbf{d}}$-action on the doubled quiver considered in \cite{CB}.
\end{enumerate}

The momentum map $\mu:R_\mathbf{d}(Q)\rightarrow\bigoplus_{i\in Q_0}{\rm Herm}(V_i)$ is equivariant for the $U_\mathbf{d}$-action, and the $U_\mathbf{d}$-orbits in $\bigoplus_i{\rm Herm}(V_i)$ are parametrized by the ordered tuples of real eigenvalues $(a_{i,1}\leq\ldots\leq a_{i,d_i})_{i\in Q_0}$ 
of the Hermitian operators. Thus, to determine the image of $\mu$, it suffices to describe the possible ordered tuples $(a_{i,k})$ of eigenvalues. A trivial restriction is
$$\sum_{i\in Q_0}\sum_{k=1}^{d_i}a_{i,k}=0$$
since the sum of traces of operators is zero on the image of $\mu$ by the previous formula.\\[2ex]
We need some combinatorial notation to formulate our main result. We consider tuples $\mathbf{K}=(K_i)_{i\in Q_0}$ of finite sets, in particular we define $$[\mathbf{d}]=(\{1,\ldots,d_i\})_i$$ for $\mathbf{d}\in\mathbb{N}Q_0$. We define the weight of $\mathbf{K}$  as the dimension vector $$|\mathbf{K}|=(|K_i|)_i\in\mathbb{N}Q_0.$$
If $\mathbf{L}\subset[|\mathbf{K}|]$, writing $K_i=\{k_{i,1}<\ldots<k_{i,|K_i|}\}$, we define the height $${\rm ht}(\mathbf{K},\mathbf{L})=\sum_{i\in Q_0}\sum_{j\in L_i}(k_{i,j}-j).$$

Our main result is:

\begin{theorem}\label{mainintro} Let $(B_i)_{i\in Q_0}$ be a tuple of Hermitian operators in $\bigoplus_i{\rm Herm}(V_i)$, and let $a_{i,1}\leq\ldots\leq a_{i,d_i}$ be the eigenvalues of $B_i$, for all $i\in Q_0$. Then $(B_i)_i$ belongs to the image of $\mu_Q$ if and only if $$\sum_{i\in Q_0}\sum_k a_{i,k}=0\mbox{ and }a_\mathbf{K}=\sum_{i\in Q_0}\sum_{k\in K_i}a_{i,k}\geq 0\mbox{ for all }\mathbf{K}\in \hat{S}(\mathbf{d}),$$ where the sets $\hat{S}(\mathbf{d})$ are defined inductively as follows:\\[1ex]
a tuple of finite sets $\mathbf{K}$ belongs to $\hat{S}(\mathbf{d})$ if and only if
$${\rm ht}(\mathbf{K},\mathbf{L})\geq\sum_{\alpha:i\rightarrow j}|L_i|(d_j-|K_j|)\mbox{ for all }\mathbf{L}\in \hat{S}(|\mathbf{K}|).$$
\end{theorem}

\section{General subrepresentations}\label{general}

Let $\mathbf{e}\leq\mathbf{d}$ be a dimension vector with entries $e_i\leq d_i$ for all $i\in Q_0$. We say that $\mathbf{e}$ is a general subrepresentation vector, written $\mathbf{e}\hookrightarrow\mathbf{d}$, if there exists a Zariski-dense subset $U$ of $R_\mathbf{d}(Q)$ such that every representation $V=((V_i)_i,(f_\alpha)_\alpha)$ for $(f_\alpha)_\alpha\in U$ admits a subrepresentation of dimension vector $\mathbf{e}$.\\[1ex]
This can be interpreted geometrically as follows: inside the product $$R_\mathbf{d}(Q)\times\prod_{i\in Q_0}{\rm Gr}_{e_i}(V_i)$$ of the representation space with a product of Grassmannians, we consider the Zariski-closed subset ${\rm Gr}_\mathbf{e}^Q(\mathbf{d})$ of pairs $((f_\alpha)_\alpha,(U_i)_i)$ such that $$f_\alpha(U_i)\subset U_j\mbox{ for all }\alpha:i\rightarrow j$$ (called the universal quiver Grassmannian). The projection $${\rm Gr}_\mathbf{e}^Q(\mathbf{d})\rightarrow \prod_i{\rm Gr}_{e_i}(V_i)$$ is a homogeneous vector bundle, thus ${\rm Gr}_\mathbf{e}^Q(\mathbf{d})$ is a smooth irreducible variety. The projection $\pi:{\rm Gr}_\mathbf{e}^Q(\mathbf{d})\rightarrow R_\mathbf{d}(Q)$ is proper, thus has Zariski-closed image. By definition $\mathbf{e}\hookrightarrow\mathbf{d}$ is equivalent to $\pi$ having dense image, thus, by properness, to $\pi$ being surjective.\\[1ex]
For two arbitrary dimension vectors $\mathbf{d}$, $\mathbf{e}$, we denote by ${\rm ext}(\mathbf{d},\mathbf{e})$ the minimal value of $\dim{\rm Ext}^1(V,W)$ for representations $V$ and $W$ of dimension vectors $\mathbf{d}$ and $\mathbf{e}$, respectively. This minimal value is in fact assumed on a Zariski-open subset of $R_\mathbf{d}(Q)\times R_\mathbf{e}(Q)$.\\[1ex]
By a theorem of Schofield \cite[Theorem 3.3]{Sch}, we have:
$$\mathbf{e}\hookrightarrow\mathbf{d}\mbox{ if and only if }{\rm ext}(\mathbf{e},\mathbf{d}-\mathbf{e})=0.$$

By another theorem of Schofield \cite[Theorem 5.4]{Sch}, we have
$${\rm ext}(\mathbf{d},\mathbf{e})=-\max(\langle\mathbf{d'},\mathbf{e}\rangle\, |\, \mathbf{d'}\hookrightarrow\mathbf{d}),$$
where we use the Euler form of the quiver defined by
$$\langle\mathbf{d},\mathbf{e}\rangle=\sum_{i\in Q_0}d_ie_i-\sum_{\alpha:i\rightarrow j}d_ie_j.$$

Thus, whether $\mathbf{e}\hookrightarrow\mathbf{d}$ holds can be determined recursively: denote by $S_Q(\mathbf{d})$ the set of all $\mathbf{e}$ such that $\mathbf{e}\hookrightarrow\mathbf{d}$; then
$$S_Q(\mathbf{d})=\{\mathbf{e}\leq\mathbf{d}\, |\, \langle \mathbf{e'},\mathbf{d}-\mathbf{e}\rangle\geq 0\mbox{ for all }\mathbf{e'}\in S_Q(\mathbf{e})\}.$$

We apply these notions to determine the scalar part of the momentum map image:\\[1ex]
Call an element of ${\bigoplus_i}{\rm Hilb}(V_i)$ scalar if it is a tuple $(\Theta_i\cdot{\rm id}_{V_i})_i$ of scalar operators for $\Theta_i\in\mathbb{R}$. We define $\Theta(\mathbf{e})=\sum_i\Theta_ie_i$. By a result of Chindris \cite[Proposition 1.3]{Ch} generalizing \cite[Proposition 6.5]{Kin}, the following are equivalent:
\begin{enumerate}
\item The tuple $(\Theta_i\cdot{\rm id}_{V_i})_i$ belongs to the image of $\mu$,
\item there exists a $\Theta$-semistable representation of dimension vector $\mathbf{d}$ in the sense of \cite{Kin},
\item we have $\Theta(\mathbf{d})=0$ and $\Theta(\mathbf{e})\geq 0$ for all $\mathbf{e}\hookrightarrow\mathbf{d}$.
\end{enumerate}

Using the above notion, we thus find:

\begin{proposition}\label{king} The scalar part of the momentum map image consists of all tuples $(\Theta_i\cdot{\rm id}_{V_i})_i$ such that $\Theta(\mathbf{d})=0$ and $\Theta(\mathbf{e})\geq 0$ for all $\mathbf{e}\in S_Q(\mathbf{d})$.
\end{proposition}

\section{From $Q$ to the leg-extended quiver $Q_\mathbf{d}$}\label{leg}

Given $Q$ and $\mathbf{d}$ as before, we define the leg-extended quiver $Q_\mathbf{d}$ as the quiver with vertices $(i,k)$ for $i\in Q_0$ and $k=1,\ldots,d_i$, and arrows $$\alpha:(i,d_i)\rightarrow (j,d_j)\mbox{ for all }\alpha:i\rightarrow j\mbox{ in }Q,$$
as well as arrows $$\beta_{i,k}:(i,k)\rightarrow(i,k+1)\mbox{ for all }i\in Q_0\mbox{ and }k=1,\ldots,d_i-1.$$
The quiver $Q_\mathbf{d}$ is thus obtained from $Q$ by attaching a leg of length $d_i-1$ at every vertex $i\in Q_0$; 
this quiver plays a prominent role in \cite{HLV}. We also define a dimension vector $\mathbf{\hat{d}}$ for $Q_\mathbf{d}$ by $$\mathbf{\hat{d}}_{i,k}=k.$$ For a dimension vector $\mathbf{{e}}$ of $Q_\mathbf{d}$, we denote by $|\mathbf{{e}}|$ the dimension vector of $Q$ given by $$|\mathbf{{e}}|_i={e}_{i,d_i}.$$

A representation $\hat{V}$ of $Q_\mathbf{\mathbf{d}}$ of dimension vector $\mathbf{{e}}\leq\mathbf{\hat{d}}$ is thus given by a representation $V$ of $Q$ of dimension vector $\mathbf{d}$, together with a chain of maps $$V_{i,1}\rightarrow V_{i,2}\rightarrow\ldots\rightarrow V_{i,d_i}=V_i$$ for every $i\in Q_0$.\\[1ex]
Given a Hermitian $n\times n$-matrix $B$, nonpositive real numbers $\Theta_1,\ldots,\Theta_n$, and complex $(i+1)\times i$-matrices $A_i$, consider the following system of equations (where $E_k$ denotes a $k\times k$ identity matrix):
$$-A_1^*A_1=\Theta_1\cdot E_1,$$
$$A_iA_i^*-A_{i+1}^*A_{i+1}=\Theta_{i+1}\cdot E_{i+1}\mbox{ for }i=1,\ldots,n-2,$$
$$A_{n-1}A_{n-1}^*+B=\Theta_n\cdot E_n.$$
The following key lemma is contained in \cite{CG}:

\begin{lemma} In a solution to the above system of equations, the matrix $B$ has eigenvalues $a_i=\Theta_i+\ldots+\Theta_n$ for $i=1,\ldots,n$. Conversely, given a Hermitian matrix $B$ with these eigenvalues, there exist matrices $A_1,\ldots,A_{n-1}$ solving the above system of equations.
\end{lemma}

\proof To prove the first claim, we use the following identity holding for arbitrary pairs of a complex $m\times n$-matrix $A$ and a complex $n\times m$-matrix $B$:
$$\det(AB-\lambda\cdot E_m)=(-\lambda)^{m-n}\det(BA-\lambda\cdot E_n)$$
(which is easily verified by assuming, via base change, the matrix $A$ to be of block form $$\left[\begin{array}{ll}E_r&0\\ 0&0\end{array}\right]).$$
Using this, we can calculate:
$$\det(B-\lambda\cdot E_n)=\det(\Theta_n\cdot E_n-A_{n-1}A_{n-1}^*-\lambda\cdot E_n)=$$
$$=\det(-A_{n-1}A_{n-1}^*+(\Theta_n-\lambda)\cdot E_n)=$$
$$=(-1)^n\det(A_{n-1}A_{n-1}^*-(\Theta_n-\lambda)\cdot E_n)=$$
$$=(-1)^n\det(A_{n-1}^*A_{n-1}-(\Theta_n-\lambda)\cdot E_{n-1})\cdot(\lambda-\Theta_n)=$$
$$=(-1)^n(\lambda-\Theta_n)\det(A_{n-2}A_{n-2}^*-\Theta_{n-1}\cdot E_{n-1}-(\Theta_n-\lambda)\cdot E_{n-1})=$$
$$=(-1)^n(\lambda-\Theta_n)\det(A_{n-2}A_{n-2}^*-(\Theta_{n-1}+\Theta_n-\lambda)\cdot E_{n-1})=$$
$$=(-1)^n(\lambda-\Theta_n)(\lambda-\Theta_{n-1}-\Theta_n)\det(A_{n-2}^*A_{n-2}-(\Theta_{n-1}+\Theta_n-\lambda)\cdot E_{n-2})=$$
$$\ldots$$
$$=(-1)^n(\lambda-\Theta_n)\cdot\ldots\cdot(\lambda-\Theta_1-\ldots-\Theta_n).$$
Thus $B$ has eigenvalues $a_1,\ldots,a_n$, as claimed.\\[1ex]
Conversely, assume that $B$ has eigenvalues $a_1\leq\ldots\leq a_n$. By equivariance of the above system of equations for the natural action of the product of unitary groups $U_1\times\ldots\times U_n$, we can assume $B$ to be the diagonal matrix ${\rm diag}(a_1,\ldots,a_n)$. For all $j\leq i$, we define $$c_i^j=\sqrt{a_{i+1}-a_j},$$ and we define $A_i$ as the $(i+1)\times i$-matrix with entries $$(A_i)_{jj}=c_i^j\mbox{ for }j\leq i,\mbox{ and }(A_i)_{jk}=0\mbox{ otherwise}.$$

The above system of equations is easily verified, proving the lemma.\\[1ex]
Using this lemma, we can identify the momentum map image for $Q$ for dimension vector $\mathbf{d}$ with the scalar part of the momentum map image for $Q_\mathbf{d}$ for the dimension vector $\mathbf{\hat{d}}$ as follows:

\begin{proposition}\label{extend} Let $(B_i)_{i\in Q_0}$ be a tuple of Hermitian operators in $\bigoplus_i{\rm Herm}(V_i)$, and let $a_{i,1}\leq\ldots\leq a_{i,d_i}$ be the eigenvalues of $B_i$, for all $i\in Q_0$. Set $$\Theta_{i,k}=a_{i,k}-a_{i,k+1}\mbox{ for all }i\in Q_0\mbox{ and }k<d_i,\mbox{ and }\Theta_{i,d_i}=a_{i,d_i}\mbox{ for all }i\in Q_0.$$ Then $(B_i)_i$ belongs to the image of $\mu_Q$ if and only if $(\Theta_{i,k}\cdot {\rm Id}_{V_{i,k}})_{i,k}$ belongs to the scalar part of the image of $\mu_{Q_\mathbf{d}}$.
\end{proposition}

\proof A representation $\hat{V}$ of $Q_\mathbf{d}$ of dimension vector $\mathbf{\hat{d}}$ on vector spaces $V_{i,k}$ is given by a representation $V$ of $Q$ of dimension vector $\mathbf{d}$ on the vector spaces $V_i=V_{i,d_i}$, given by linear maps $(f_\alpha:V_i\rightarrow V_j)_{\alpha:i\rightarrow j}$, together with maps $g_{i,k}:V_{i,k}\rightarrow V_{i,k+1}$ for all $i\in Q_0$ and $k<d_i$. The group ${\rm GL}_\mathbf{\hat{d}}$ acts by base change in all vector spaces $V_{i,k}$. The momentum map for $Q_\mathbf{d}$ is then given as follows:
$$\mu_{Q_\mathbf{d}}(\hat{V})_{i,1}=-g_{i,1}^*\circ g_{i,1},$$
$$\mu_{Q_\mathbf{d}}(\hat{V})_{i,k}=g_{i,k-1}\circ g_{i,k-1}^*-g_{i,k}^*\circ g_{i,k}\mbox{ for }k=2,\ldots,d_i-1,$$
$$\mu_{Q_\mathbf{d}}(\hat{V})_{i,d_i}=\mu_Q(V)_i+g_{i,d_i-1}\circ g_{i,d_i-1}^*$$
for all $i\in Q_0$.\\[1ex]
If $(\Theta_{i,k}\cdot {\rm id}_{V_{i,k}})_{i,k}$ belongs to the scalar part of the image of $\mu_{Q_\mathbf{d}}$, namely as the image of a representation $\hat{V}$, then, by the previous lemma, each $\mu_Q(V)_i$ has eigenvalues $a_{i,1},\ldots,a_{i,d_i}$. Thus $\mu_Q(V)$ is conjugate to $(B_i)_i$ under $U_\mathbf{d}$. By $U_\mathbf{d}$-equivariance of $\mu_Q$, the tuple $(B_i)_i$ also belongs to the image of $\mu_Q$. Conversely, assume that $(B_i)_i$ belongs to the image of $\mu_Q$, namely as the image of a representation $V$. Again by the previous lemma, we find linear maps $g_{i,k}$, thus completing $V$ to a representation $\hat{V}$ of $Q_\mathbf{d}$, such that $\mu_{Q_\mathbf{d}}(\hat{V})$ is given by $(\Theta_{i,k}\cdot{\rm id}_{V_{i,k}})_{i,k}$, as claimed.\\[1ex]
We rephrase this result in terms of the following diagram:
\begin{center}
\begin{tikzcd}
R_\mathbf{d}(Q) \arrow[r, "\mu_{Q}"] 
& \bigoplus_{i\in Q_0}{\rm Herm}(V_i) \arrow[dr, leftarrow,  "S" ] 
\\ R_\mathbf{\hat{d}}(Q_\mathbf{d}) \arrow[r]
& \bigoplus_{i\in Q_0}\bigoplus_{k\leq d_i}{\rm Herm}(V_{i,k}) \arrow[r, hookleftarrow]
& \bigoplus_{i\in Q_0}\bigoplus_{k\leq d_i}\mathbb{R}\cdot{\rm Id}_{V_{i,k}},
\end{tikzcd}
\end{center}

where $S$ maps a tuple $(\Theta_{i,k}\cdot{\rm Id}_{V_{i,k}})$ to the tuple of diagonal matrices $B_i$ with eigenvalues $a_{i,k}$, respectively. Then the previous proposition can be rephrased as
$$\mu_Q(R_\mathbf{d}(Q))=U_\mathbf{d}\cdot S\left(\mu_{Q_\mathbf{d}}(R_\mathbf{\hat{d}}(Q_\mathbf{d}))\cap\bigoplus_{i\in Q_0}\bigoplus_{k\leq d_i}\mathbb{R}\cdot{\rm Id}_{V_{i,k}}\right).$$

\section{Reductions}\label{reductions}

We can now give a first explicit description of the image of the momentum map using the result of the previous section.\\[1ex] 
As before, let $(B_i)_{i\in Q_0}$ be a tuple of Hermitian operators in $\bigoplus_i{\rm Herm}(V_i)$, and let $a_{i,1}\leq\ldots\leq a_{i,d_i}$ be the eigenvalues of $B_i$, for all $i\in Q_0$. Set $$\Theta_{i,k}=a_{i,k}-a_{i,k+1}\mbox{ for all }i\in Q_0\mbox{ and }k<d_i,\mbox{ and }\Theta_{i,d_i}=a_{d_i}\mbox{ for all }i\in Q_0,$$ and denote $\hat{\Theta}=(\Theta_{i,k})_{i,k}$. Combining Proposition \ref{king} and Proposition \ref{extend}, we find immediately:

\begin{corollary}\label{cor51}
The tuple $(B_i)_i$ belongs to the image of $\mu_Q$ if and only if $\hat{\Theta}(\mathbf{\hat{d}})=0$ and $\hat{\Theta}(\mathbf{{e}})\geq 0$ for all $\mathbf{{e}}\in S_{Q_\mathbf{d}}(\mathbf{\hat{d}})$, where the sets $S_{Q_\mathbf{d}}(\mathbf{{e}})$ are defined inductively by $$S_{Q_\mathbf{d}}(\mathbf{{e}})=\{\mathbf{{e}'}\leq\mathbf{{e}}\, |\, \langle \mathbf{{e}''},\mathbf{{e}}-\mathbf{{e}'}\rangle\geq 0\mbox{ for all }\mathbf{{e}''}\in S_{Q_\mathbf{d}}(\mathbf{{e}'})\}.$$
\end{corollary}

We will now perform several reductions to make this criterion more effective.

\subsection{Reduction to dimension vectors of flag type}

If all the maps representing the arrows $\beta_{i,k}$, in a representation $\hat{V}$ of $Q_\mathbf{d}$ of dimension vector $\mathbf{e}$ 
are injective, we call $\hat{V}$ {\it of flag type}, since the images $F_i^k$ of the $V_{i,k}$ in $V_i$ induce flags $F_i^*$ in every $V_i$. Moreover, if $\mathbf{{e}}=\mathbf{\hat{d}}$, the flags $F_i^*$ are complete, that is, $\dim F_i^k=k$ for all $i\in Q_0$ and $k\leq d_i$.\\[1ex]
We define a dimension vector $\mathbf{{e}}\leq\mathbf{\hat{d}}$ to be {\it of flag type} if $${e}_{i,k}-{e}_{i,k-1}\in\{0,1\}\mbox{ for all }i\in Q_0\mbox{ and }k=1,\ldots,d_i$$ (where we formally set ${e}_{i,0}=0$).

\begin{lemma}\label{flagtype} If $\mathbf{{e}}\leq\mathbf{\hat{d}}$ is of flag type and $\mathbf{{e}'}\hookrightarrow\mathbf{{e}}$, there exists a unique dimension vector $\mathbf{{e}''}$ of flag type such that $|\mathbf{e}''|=|\mathbf{e}'|$, $\mathbf{{e}'}\leq\mathbf{{e}''}\hookrightarrow\mathbf{{e}}$ and $$\langle \mathbf{{e}'},\mathbf{\hat{d}}-\mathbf{{e}}\rangle\geq \langle \mathbf{{e}''},\mathbf{\hat{d}}-\mathbf{{e}}\rangle.$$
\end{lemma}

\proof As discussed above, the condition $\mathbf{{e}'}\hookrightarrow\mathbf{{e}}$ is equivalent to the map $$\hat{\pi}:{\rm Gr}_\mathbf{{e}'}^{Q_\mathbf{d}}(\mathbf{{e}})\rightarrow R_\mathbf{\hat{e}}(Q_\mathbf{d})$$ being surjective. Restricting to the Zariski-open subsets of representations of flag type in both varieties, we have a surjective map $$\hat{\pi}:{\rm Gr}_\mathbf{{e}'}^{Q_\mathbf{d}}(\mathbf{{e}})^{\rm flag}\rightarrow R_\mathbf{{e}}(Q_\mathbf{d})^{\rm flag}.$$

An element of ${\rm Gr}_\mathbf{{e}'}^{Q_\mathbf{d}}(\mathbf{{e}})^{\rm flag}$ consists of a representation $N$ of $Q$ of dimension vector $|\mathbf{{e}}|$, flags $F_i^*$ in $N_i$ for all $i\in Q_0$, and subspaces $U_{i,k}\subset F_i^k$ of dimension ${e}'_{i,k}$ for all $i\in Q_0$ and $k\leq d_i$ forming flags $$U_{i,1}\subset\ldots\subset U_{i,d_i}=:U_i\subset N_i$$ for all $i\in Q_0$.\\[1ex]
We stratify ${\rm Gr}_\mathbf{{e}'}^{Q_\mathbf{d}}(\mathbf{{e}})^{\rm flag}$ into Zariski-locally closed strata $S_\mathbf{{f}}$ defined by the conditions $$\dim U_i\cap F_i^k={f}_{i,k}$$ for all $i\in Q_0$ and $k\leq d_i$. Suppose that a stratum $S_\mathbf{{f}}$ 
is non-empty. Noting that $ U_i\cap F_i^{d_i}=U_i$ is of dimension ${e}'_{i,d_i}$, we find $|\mathbf{{f}}|=|\mathbf{{e}'}|$. Furthermore, we have
$$(U_i\cap F_i^k)/(U_i\cap F_i^{k-1})\hookrightarrow F_i^k/F_i^{k-1}.$$
The right hand side being at most one-dimensional by $\mathbf{{e}}$ being of flag type, we see that $\mathbf{{f}}$ is of flag type, too.\\[1ex]
The variety ${\rm Gr}_\mathbf{{e}'}^{Q_\mathbf{d}}(\mathbf{{e}})$ being irreducible, precisely one of the strata $S_\mathbf{{f}}$ must be dense; define $\mathbf{{e}''}$ as the dimension vector $\mathbf{{f}}$ for which this holds. Then $\hat{\pi}(S_\mathbf{{e}''})$ is dense in $R_\mathbf{{e}}(Q_\mathbf{d})^{\rm flag}$, which means that for almost all $(N,(F_i^*)_i)$ in $R_\mathbf{{e}}(Q_\mathbf{d})^{\rm flag}$, there exist subspaces $U'_{i,k}$ (namely the $U_i\cap F_i^k$) of dimension ${e}''_{i,k}$ of the $F_i^k$ forming flags inside each $N_i$. But this precisely means that $\mathbf{{e}''}\hookrightarrow\mathbf{{e}}$.\\[1ex]
Denoting $\mathbf{{g}}=\mathbf{{e}'}-\mathbf{{e}''}\geq 0$, it remains to prove that $\langle\mathbf{{g}},\mathbf{{d}}-\mathbf{{e}}\rangle\leq 0$. Since $|\mathbf{{g}}|=0$, we can easily calculate
$$\langle\mathbf{{g}},\mathbf{\hat{d}}-\mathbf{{e}}\rangle=\sum_{i\in Q_0}\sum_{k<d_i}{g}_{i,k}({e}_{i,k+1}-{e}_{i,k}-1),$$
which is nonpositive since $\mathbf{{e}}$ is of flag type, finishing the proof.\\[1ex]
We claim that, in the description of the momentum map image, we can reduce to the consideration of dimension vectors of flag type. More precisely, we have:

\begin{proposition}\label{almost} Define sets $S_{Q_\mathbf{d}}(\mathbf{{e}})^{\rm flag}$ inductively by
$$S_{Q_\mathbf{d}}(\mathbf{{e}})^{\rm flag}=\{\mathbf{{e}'}\leq\mathbf{{e}}\mbox{ of flag type }\, |\, \langle \mathbf{{e}''},\mathbf{{e}}-\mathbf{{e}'}\rangle\geq 0\mbox{ for all }\mathbf{{e}''}\in S_{Q_\mathbf{d}}(\mathbf{{e}'})^{\rm flag}\}.$$
Then a tuple $(B_i)_i$ belongs to the image of $\mu_Q$ if and only if for the associated $\hat{\Theta}$ as above, we have $\hat{\Theta}(\mathbf{\hat{d}})=0$ and $\hat{\Theta}(\mathbf{{e}})\geq 0$ for all $\mathbf{{e}}\in S_{Q_\mathbf{d}}(\mathbf{\hat{d}})^{\rm flag}$.
\end{proposition}

\proof Assume $\hat{\Theta}$ is given as above, thus in particular $\Theta_{i,k}\leq 0$ for all $i\in Q_0$ and all $k<d_i$. In light of the previous statement, we first have to prove that $\hat{\Theta}(\mathbf{{e}})\geq 0$ holds for all $\mathbf{{e}}\in S_{Q_\mathbf{d}}(\mathbf{\hat{d}})$ if and only if it holds for all $\mathbf{{e}}\in S_{Q_\mathbf{d}}(\mathbf{\hat{d}})^{\rm flag}$. If $\mathbf{{e}}\in S_{Q_\mathbf{d}}(\mathbf{\hat{d}})$ is not of flag type, we find by Lemma \ref{flagtype} a dimension vector $\mathbf{{e}'}\geq\mathbf{{e}}$ of flag type such that $|\mathbf{{e}'}|=|\mathbf{{e}}|$. This implies that $\hat{\Theta}(\mathbf{{e}})\geq\hat{\Theta}(\mathbf{{e}'})$, thus the condition $\hat{\Theta}(\mathbf{{e}})\geq 0$ is superfluous given the conditions $\hat{\Theta}(\mathbf{{e}'})\geq 0$ for all $\mathbf{{e}'}\in S_{Q_\mathbf{d}}(\mathbf{\hat{d}})^{\rm flag}$.\\[1ex]
Now we have to prove that $\mathbf{{e}'}\hookrightarrow\mathbf{{e}}\leq\mathbf{\hat{d}}$ holds for dimension vectors of flag type if and only if $\langle\mathbf{{e}''},\mathbf{{e}}-\mathbf{{e}'}\rangle\geq 0$ for all $\mathbf{{e}''}$ of flag type such that $\mathbf{{e}''}\hookrightarrow\mathbf{{e}'}$. If such a dimension vector $\mathbf{{e}''}$ is not of flag type, again by Lemma \ref{flagtype} we find a dimension vector $\mathbf{{e}'''}$ of flag type such that $\mathbf{{e}''}\leq\mathbf{{e}'''}$, and thus $$\langle \mathbf{{e}''},\mathbf{{e}}-\mathbf{{e}'}\rangle\geq \langle\mathbf{{e}'''},\mathbf{{e}}-\mathbf{{e}'}\rangle\geq 0,$$
finishing the proof.

\subsection{Reduction to smaller leg-extended quivers}

Given dimension vectors of flag type $\mathbf{e}$ for $Q_\mathbf{d}$ and $\mathbf{f}$ for $Q_{|\mathbf{e}|}$, we define $\exp_\mathbf{e}\mathbf{f}$, the $\mathbf{e}$-expansion of $\mathbf{f}$, as the dimension vector for $Q_\mathbf{d}$ given by
$$(\exp_\mathbf{e}\mathbf{f})_{i,k}=f_{i,j}\mbox{ if }e_{i,k}=j.$$

\begin{lemma} Let $\mathbf{e}\leq\mathbf{\hat{d}}$ be a dimension vector of flag type for $Q_\mathbf{d}$. Then, for a dimension vector $\mathbf{e'}\leq\mathbf{e}$ of flag type, the following are equivalent:
\begin{enumerate}
\item $\mathbf{e'}$ fulfills $\mathbf{e'}\hookrightarrow\mathbf{e}$,
\item $\mathbf{e'}\leq\exp_\mathbf{e}\mathbf{f}$ and $|\mathbf{e'}|=|\mathbf{f}|$ for a dimension vector $\mathbf{f}\leq\widehat{|\mathbf{e}|}$ for $Q_{|\mathbf{e}|}$.
\end{enumerate}
\end{lemma}

\proof The first condition means that, for all representations $W$ in $R_{|\mathbf{e}|}(Q)$ and all $Q_0$-tuples of flags $$F_{i,1}\subset\ldots\subset F_{i,d_i}=W_i$$ of dimension vector $\mathbf{e}$, there exists a subrepresentation $U$ of $W$ of dimension vector
$|\mathbf{e'}|$ such that $$\dim U_i\cap F_{i,k}\geq e'_{i,k}$$ for all $i\in Q_0$ and $k=1,\ldots,d_i$. Then
$$(U_i\cap F_{i,k})/(U_i\cap F_{i,k-1})\hookrightarrow F_{i,k}/F_{i,k-1},$$
thus the dimensions $\dim U_i\cap F_{i,k}$ form a dimension vector $\mathbf{e''}$ of flag type for $Q_\mathbf{d}$, which is constant along the leg segments where $\mathbf{e}$ is constant. This precisely means that $\mathbf{e''}=\exp_\mathbf{e}\mathbf{f}$ for a dimension vector of flag type $\mathbf{f}$ for $Q_{|\mathbf{e}|}$. Contracting the flags $(F_{i,k})_{i,k}$ to flags of dimension vector $\mathbf{f}$ for $Q_{|\mathbf{e}|}$, we see that $\mathbf{f}\hookrightarrow\widehat{|\mathbf{e}|}$ for $Q_{|\mathbf{e}|}$.\\[1ex]
Conversely, the second condition means that, for all representations $W$ in $R_{|\mathbf{e}|}$ and all $Q_0$ tuples of flags $\bar{F}_{i,k}$ of dimension vector $\widehat{|\mathbf{e}|}$ for $Q_{|\mathbf{e}|}$, there exists a subrepresentation $U$ of dimension vector $|\mathbf{e'}|$ such that $$\dim U_i\cap \bar{F}_{i,k}\geq f_{i,k}$$ for all $i$ and $k$. Furthermore, we have $\mathbf{e'}\leq\exp_\mathbf{e}\mathbf{f}$ and $|\mathbf{e'}|=|\mathbf{f}|$. Expanding the flags $\bar{F}_i$ to flags $F_i$ for $Q_\mathbf{d}$ which are constant along the leg segments where $\mathbf{e}$ is constant, we thus find $$(\dim U_i\cap F_{i,k})_{i,k}\geq\exp_\mathbf{e}\mathbf{f}\geq\mathbf{e'},$$ and thus $\mathbf{e'}\hookrightarrow\mathbf{e}$, finishing the proof.\\[1ex]
We now define sets $\check{S}(\mathbf{d})$ of dimension vectors for $Q_\mathbf{d}$ recursively by
$$\check{S}(\mathbf{d})=\{\mathbf{e}\leq\mathbf{\hat{d}}\mbox{ of flag type }\, |\, \langle \exp_\mathbf{e}\mathbf{f},\mathbf{\hat{d}}-\mathbf{e}\rangle\geq 0\mbox{ for all }\mathbf{f}\in\check{S}(|\mathbf{e}|)\}.$$

\begin{proposition}\label{propq} We have $\check{S}(\mathbf{d})=S_{Q_\mathbf{d}}(\mathbf{\hat{d}})^{\rm flag}$ for all $\mathbf{d}$.
\end{proposition}

\proof The proof proceeds by induction over $\mathbf{d}$, and we prove both inclusions separately. Let $\mathbf{e}\leq\mathbf{\hat{d}}$ be of flag type such that $\langle\exp_\mathbf{e}\mathbf{f},\mathbf{\hat{d}}-\mathbf{e}\rangle\geq 0$ for all $\mathbf{f}\in\check{S}(|\mathbf{e}|)$. We can assume $\mathbf{e}\not=\mathbf{\hat{d}}$, thus, by the inductive assumption, we have $$\check{S}(|\mathbf{e}|)=S_{Q_{|\mathbf{e}|}}(\widehat{|\mathbf{e}|})^{\rm flag},$$ and thus $\langle\exp_\mathbf{e}\mathbf{f},\mathbf{\hat{d}}-\mathbf{e}\rangle\geq 0$ for all $\mathbf{f}\hookrightarrow\widehat{|\mathbf{e}|}$ of flag type. We have to show that $\mathbf{e}\hookrightarrow\mathbf{\hat{d}}$, or, equivalently, $\langle\mathbf{e}',\mathbf{\hat{d}}-\mathbf{e}\rangle\geq 0$ for all $\mathbf{e}'\hookrightarrow\mathbf{e}$. Let such a $\mathbf{e'}$ be given. By the previous lemma, we then find a dimension vector $\mathbf{f}\hookrightarrow\widehat{|\mathbf{e}|}$ of flag type such that $\mathbf{e'}\leq\exp_\mathbf{e}\mathbf{f}$ and $|\mathbf{e'}|=|\exp_\mathbf{e}\mathbf{f}|=|\mathbf{f}|$. Thus $$0\leq\langle\exp_\mathbf{e}\mathbf{f},\mathbf{\hat{d}}-\mathbf{e}\rangle\leq\langle\mathbf{e'},\mathbf{\hat{d}}-\mathbf{e}\rangle$$ as desired.\\[1ex]
Conversely, let $\mathbf{e}\hookrightarrow\mathbf{\hat{d}}$ be of flag type, and let $\mathbf{f}$ be in $\check{S}(\widehat{|\mathbf{e}|})$. We have to prove that $\langle\exp_\mathbf{e}\mathbf{f},\mathbf{\hat{d}}-\mathbf{e}\rangle\geq 0$. We can assume $\mathbf{e}\not=\mathbf{\hat{d}}$, thus, by the inductive assumption, we have $\mathbf{f}\in S_{Q_{|\mathbf{e}|}}(\widehat{|\mathbf{e}|})$, that is, $\mathbf{f}\hookrightarrow\widehat{|\mathbf{e}|}$. Then $\exp_\mathbf{e}\mathbf{f}\hookrightarrow\mathbf{e}$, and thus $\langle\exp_\mathbf{e}\mathbf{f},\mathbf{\hat{d}}-\mathbf{e}\rangle\geq 0$ as desired, finishing the proof.

\subsection{Reduction to tuples of sets}

For a tuple $$\mathbf{K}=(K_i\subset \{1,\ldots, d_i\})_{i\in Q_0}$$ of finite sets, we define a dimension vector $\mathbf{d}(\mathbf{K})$ for $Q_\mathbf{d}$ by $$d(\mathbf{K})_{i,k}=|K_i\cap\{1,\ldots,k\}|$$ for all $i\in Q_0$, $k=1,\ldots,d_i$. In particular, setting $$[\mathbf{d}]=(\{1,\ldots,d_i\})_i,$$ we have $$\mathbf{d}([\mathbf{d}])=\hat{\mathbf{d}}.$$

Such dimension vectors $\mathbf{d}(\mathbf{K})$ are obviously of flag type, and, conversely, all dimension vectors $\mathbf{e}$ for $Q_\mathbf{d}$ which are of flag type arise in this way; namely, defining $K_i$ as the set of indices $k$ where $e_{i,k}=e_{i,k-1}+1$ for all $i\in Q_0$, we obviously find $\mathbf{e}=\mathbf{d}(\mathbf{K})$.\\[1ex]
We can now reformulate Proposition \ref{propq} in this terminology. First of all, it is easy to verify, by telescoping the summation, that
$$\hat{\Theta}(\mathbf{d}(\mathbf{K}))=\sum_{i\in Q_0}\sum_{k\in K_i}a_{i,k}.$$

As before, we define the weight of $\mathbf{K}$ as the dimension vector $$|\mathbf{K}|=(|K_i|)_i\in\mathbb{N}Q_0.$$ If $\mathbf{L}\subset[|\mathbf{K}|]$, we write $$K_i=\{k_{i,1}<\ldots<k_{i,|K_i|}\}$$ and define $$\mathbf{K}_\mathbf{L}=(\{k_{i,j}\, |\, j\in L_i\})_i.$$ We claim that $$\exp_{\mathbf{d}(\mathbf{K})}\mathbf{d}(\mathbf{L})=\mathbf{d}(\mathbf{K}_\mathbf{L}):$$

Indeed, we have
$$\exp_{\mathbf{d}(\mathbf{K})}\mathbf{d}(\mathbf{L})_{i,k}=\mathbf{d}(\mathbf{L})_{i,j}\mbox{ if }\mathbf{d}(\mathbf{K})_{i,k}=j,$$ that is,
$$\exp_{\mathbf{d}(\mathbf{K})}\mathbf{d}(\mathbf{L})_{i,k}=|L_i\cap\{1,\ldots,j\}|\mbox{ if }|K_i\cap\{1,\ldots,k\}|=j.$$ Thus $\exp_{\mathbf{d}(\mathbf{K})}\mathbf{d}(\mathbf{L})_{i,k}$ equals the number of indices $j\in L_i$ such that $j\leq |K_i\cap\{1,\ldots,k\}|$, or, in other words, the number of indices $j\in L_i$ such that $k_{i,j}\leq k$. But this is precisely the cardinality of $(K_i)_{L_i}\cap\{1,\ldots,k\}$, which equals $\mathbf{d}(\mathbf{K}_\mathbf{L})_{i,k}$.\\[1ex]
Finally, we define the height $${\rm ht}(\mathbf{K},\mathbf{L})=\sum_{i\in Q_0}\sum_{j\in L_i}(k_{i,j}-j)$$ and compute $\langle \exp_\mathbf{e}\mathbf{f},\mathbf{\hat{d}}-\mathbf{e}\rangle$. First, again by telescoping the summation, we have
$$\langle\mathbf{d}(\mathbf{K'}),\mathbf{\hat{d}}-\mathbf{d}(\mathbf{K})\rangle=\langle |\mathbf{K'}|,\mathbf{d}-|\mathbf{K}|\rangle_Q-\sum_{i\in Q_0}|\{l<k\, :\, l\in K_i', k\not\in K_i\}|.$$
For $\mathbf{K'}=\mathbf{K}_\mathbf{L}$, this reduces to $${\rm ht}(\mathbf{K},\mathbf{L})-\sum_{i\rightarrow j}|L_i|(d_j-|K_j|)$$
since
$$|\{l<k\, :\, l\in(K_i)_{L_i}, k\not\in K_i\}|=\sum_{l\in(K_i)_{L_i}}|\{k>l\, :\, k\not\in K_i\}|=$$
$$=\sum_{j\in L_i}|\{k>k_{i,j}\, :\, k\not\in K_i\}|=\sum_{j\in L_i}\left(|\{k>k_{i,j}\}|-|\{k>k_{i,j}\, :\, k\in K_i\}|\right)=$$
$$=\sum_{j\in L_i}((d_i-k_{i,j})-(|K_i|-j))=|L_i|\cdot(d_i-|K_i|)-\sum_{j\in L_i}(k_{i,j}-j).$$

Combining this translation to tuples of finite sets with  Proposition \ref{almost} and Proposition \ref{propq}, we arrive precisely at the formulation of Theorem \ref{mainintro}, and this main result is proved.

\section{Examples}\label{examples}

We  first illustrate our main result with a few small examples.\\[1ex]
We start with the quiver $L_m$ consisting of a single vertex to which $m$ loops are attached. For a dimension $d\geq 0$, we are thus considering the diagonal adjoint action of ${\rm GL}_d$ on the space $M_{d\times d}^m$ of $m$-tuples of $d\times d$-matrices. Choosing the standard Hermitian form on $\mathbb{C}^d$, the momentum map $$\mu:M_{d\times d}^m\rightarrow{\rm Herm}(\mathbb{C}^d)$$ is given by $$\mu(A_1,\ldots,A_m)=\sum_{k=1}^m[A_k,A_k^*].$$ The set $\hat{S}(d)$ consists of all subsets $$K=\{k_1<\ldots<k_{|K|}\}$$ of $[d]$ such that $$\sum_{j\in L}(k_j-j)\geq m|L|(d-|K|)$$ for all $L$ in $\hat{S}(|K|)$. Since $k_j\leq d-|K|+j$ trivially, we always have $$\sum_{j\in L}(k_j-j)\leq |L|(d-|K|).$$ Together with an easy induction, this yields
$$\hat{S}(d)=\left\{\begin{array}{lcl}2^{[d]}&,&m=0,\\ \{\emptyset,\{d\},\ldots,\{2,\ldots,d\},[d]\}&,&m=1\\ \{\emptyset,[d]\}&,&m\geq 2.\end{array}\right.$$
In the case $m=1$, this yields the inequalities
$$a_j+\ldots+a_d\geq 0\mbox{ for all }j$$ for the eigenvalues $a_1<\ldots<a_d$ of a traceless Hermitian matrix, which are redundant: for every $k<j$, we have
$$(d-j+1)a_k<a_j+\ldots+a_d,$$
and thus
$$d(a_j+\ldots+a_d)=(j-1)(a_j+\ldots a_d)+(d-j+1)(a_j+\ldots+a_d)=$$
$$=(j-1)(a_j+\ldots+a_d)-(d-j+1)(a_1+\ldots+a_{j-1})>0.$$

We thus conclude with our main result that every traceless Hermitian matrix is of the form $[A,A^*]$.\\[1ex]
Our second example is the $m$-subspace quiver $S_m$ with vertices $i_1,\ldots,i_m,j$ and arrows $i_k\rightarrow j$ for $k=1,\ldots,m$, and the dimension vector ${\bf d}$ given by $d_{i_k}=1$ for all $k$ and $d_j=2$. We thus consider the action of $(\mathbb{C}^*)^m\times{\rm GL}_2$ on $(\mathbb{C}^2)^m$ via $$(t_1,\ldots,t_m,g)\cdot(v_1,\ldots,v_m)=(\frac{1}{t_1}gv_1,\ldots,\frac{1}{t_m}gv_m).$$
Choosing again the standard Hermitian form on $\mathbb{C}^2$, the momentum map $$\mu:(\mathbb{C}^2)^m\rightarrow\mathbb{R}^m\times{\rm Herm}(\mathbb{C}^2)$$ is given by $$\mu(v_1,\ldots,v_m)=(-v_1^*v_1,\ldots,-v_m^*v_m,\sum_kv_kv_k^*).$$ After a direct calculation, we see that the set $\hat{S}(\mathbf{d})$ consists of the tuples
$$(\emptyset,\ldots,\emptyset,\emptyset),(\emptyset,\ldots,\emptyset,\{1\}),(\emptyset,\ldots,\emptyset,\{2\}),$$
$$(\emptyset,\ldots,\emptyset,\{1\},\emptyset,\ldots,\emptyset,\{2\})\mbox{ with $\{1\}$ placed at any position from $1$ to $m$},$$
$$(K_1,\ldots,K_m,\{1,2\})\mbox{ for arbitrary $K_k\in\{\emptyset,\{1\})$}.$$

Our main result thus gives: 

\begin{corollary} A Hermitian $2\times 2$-matrix with eigenvalues $b_1\leq b_2$ is of the form $\sum_{k=1}^mv_kv_k^*$ for vectors $v_k$ of length $\lambda_k$, respectively, if and only if $$b_1+b_2=\sum_k\lambda_k^2,\, b_1\geq 0,\mbox{ and }b_2\geq\lambda_k^2\mbox{ for all }k.$$
\end{corollary}

We can also give a sufficient criterion for the momentum map to have maximal possible image. For a dimension vector $\mathbf{d}$, we denote by ${\rm supp}(\mathbf{d})$ the full subquiver with vertices $i$ such that $d_i>0$. We call a quiver strongly connected if there exists an oriented path between any two of its vertices.

\begin{corollary} For $Q$ and $\mathbf{d}$ as before, the following are equivalent:
\begin{enumerate}
\item We have $\hat{S}(\mathbf{d})=\{\emptyset,[\mathbf{d}]\}$.
\item There exists an irreducible representation of $Q$ of dimension vector $\mathbf{d}$.
\item Either ${\rm supp}(\mathbf{d})$ is a single oriented cycle and $d_i=1$ for all $i\in{\rm supp}(\mathbf{d})$, or ${\rm supp}(\mathbf{d})$ is strongly connected and not a single oriented cycle, and $$\sum_{\alpha:i\rightarrow j}d_j\geq d_i\leq\sum_{\alpha:j\rightarrow i}d_j\mbox{ for all }i\in Q_0.$$
\end{enumerate}
In this case, the image of the momentum map $\mu$ equals the set of Hermitian operators $(B_i)_i\in\bigoplus_{i\in Q_0}{\rm Herm}(V_i)$ with trace sum $\sum_i{\rm tr}(B_i)=0$.
\end{corollary}

\proof The equivalence of (2) and (3) is precisely the content of \cite[Theorem 4]{LBP}. Assume that (1) holds and suppose that there is no irreducible representation of $Q$ of dimension vector $\mathbf{d}$. Then there exists a proper dimension vector $0\not=\mathbf{e}\hookrightarrow\mathbf{d}$. We can thus find a nontrivial dimension vector $\mathbf{f}$ of flag type form $Q_\mathbf{d}$ such that $\mathbf{f}\hookrightarrow\mathbf{\hat{d}}$ and $|\mathbf{f}|=\mathbf{e}$. In particular, $\mathbf{f}\in S_{Q_\mathbf{d}}(\mathbf{\hat{d}})$, contradicting the assumption.\\[1ex]
Conversely, assume that there exists an irreducible representation of $Q$ of dimension vector $\mathbf{d}$. Suppose that the image of $\mu$ is constrained by a proper inequality, given by a dimension vector $\mathbf{f}$ of flag type such that $\mathbf{f}\hookrightarrow\mathbf{\hat{d}}$. Then $\mathbf{e}=|\mathbf{f}|\hookrightarrow\mathbf{d}$, thus $\mathbf{e}=0$ or $\mathbf{e}=\mathbf{d}$. In the former case, $\mathbf{f}=0$, and in the latter case, $\mathbf{f}=\mathbf{\hat{d}}$, since $\mathbf{f}$ is of flag type. This contradicts the assumption that $\mathbf{f}$ defines a proper inequality for the image of $\mu$.\\[1ex]
Finally, we can give a characterization of the momentum cone to have positive volume:

\begin{corollary} For $Q$ and $\mathbf{d}$ as before, the following are equivalent:
\begin{enumerate}
\item The image of the momentum map $\mu$ has positive volume in the space of Hermitian operators $(B_i)_i\in\bigoplus_{i\in Q_0}{\rm Herm}(V_i)$ with trace sum $\sum_i{\rm tr}(B_i)=0$.
\item There exists a representation of $Q_\mathbf{d}$ of dimension vector $\mathbf{\hat{d}}$ with trivial endomorphism ring.
\item There is no tuple of finite sets $\mathbf{K}$ such that both $\mathbf{K}$ and its complement $[\mathbf{d}]\setminus\mathbf{K}$ belong to $\hat{S}(\mathbf{d})$.
\end{enumerate}
\end{corollary}

\proof The implication from (1) to (3) is obvious, since otherwise the image of $\mu$ is contained in the proper hyperplane defined by $a_\mathbf{K}=0$. Assume that (3) holds, and suppose that there is no representation of $Q_\mathbf{d}$ of dimension vector $\mathbf{\hat{d}}$ with trivial endomorphism ring. By \cite[Theorem 2.2]{Sch}, there exists a proper decomposition $\mathbf{\hat{d}}=\mathbf{e}+\mathbf{f}$ such that a generic representation of $Q_\mathbf{d}$ of dimension vector $\mathbf{\hat{d}}$ admits a direct sum decomposition into subrepresentations of dimension vectors $\mathbf{e}$ and $\mathbf{f}$, respectively. Without loss of generality, we can assume $\mathbf{e}$ and $\mathbf{f}$ to be of flag type. Thus $\mathbf{e},\mathbf{f}\hookrightarrow\mathbf{\hat{d}}$, which translates into the existence of complementary tuples of sets $\mathbf{K},\mathbf{L}\in\hat{S}(\mathbf{d})$, a contradiction. Now assume that (2) holds. Then by \cite[Theorem 6.1]{Sch}, we have $\langle\mathbf{e},\mathbf{\hat{d}}\rangle-\langle\mathbf{\hat{d}},\mathbf{e}\rangle>0$ for all $\mathbf{e}\hookrightarrow\mathbf{\hat{d}}$. Thus $\hat{\Theta}=\langle\_,\mathbf{\hat{d}}\rangle-\langle\mathbf{\hat{d}},\_\rangle$ defines a linear form on dimension vectors for $Q_\mathbf{d}$ which fulfills $\hat{\Theta}(\mathbf{e})>0$ for all $\mathbf{e}\in S_{Q_\mathbf{d}}(\mathbf{\hat{d}})$. By our main result, this translates into the existence of an interior point in the image of $\mu$.\\[1ex]
Existence of a quiver representation with trivial endomorphism ring can again be characterized recursively using the methods of \cite{Sch}; here we derive at least a simple non-recursive sufficient criterion:

\begin{corollary} Assume that $d_i\not=0$ for all $i\in Q_0$, and that $d_i\leq\sum_{j\leftrightarrow i}d_j-1$ for all $i\in Q_0$, with strict inequality for at least on $i\in Q_0$. Then the image of the momentum map $\mu$ has positive volume.
\end{corollary}

\proof The stated inequality is equivalent to $\langle\mathbf{\hat{d}},e_{i,k}\rangle+\langle e_{i,k},\mathbf{\hat{d}}\rangle\leq 0$ for all coordinate vectors $e_{i,k}$ for $Q_\mathbf{d}$, with strict inequality for at least one of the vertices $(i,d_i)$. In the language of \cite{Kac}, this means that $\mathbf{\hat{d}}$ belongs to the so-called fundamental domain of $Q_\mathbf{d}$ and has non-tame support. In this case, \cite[Lemma 2.7]{Kac} ensures that there exists a representation with trivial endomorphism ring of $Q_\mathbf{d}$ of dimension vector $\mathbf{\hat{d}}$.

\section{Real representation spaces}\label{real}

The inductive description of the momentum map image can be easily transferred to the situation of real representation spaces of quivers using the general result \cite{OS}. We consider complex conjugation as an involution both on the group $U_\mathbf{d}$ and its Lie algebra, as well as on the space $R_\mathbf{d}(Q)$. Passing to fixed points under complex conjugation yields the following setup:\\[1ex]
Given real vector spaces $W_i$ for all vertices $i$ of the quiver $Q$, a choice of a scalar product on every vector space $W_i$ yields a scalar product on $R_\mathbf{d}(Q)$ via
$$((f_\alpha),(g_\alpha))=\sum_\alpha{\rm trace}(f_\alpha\circ g_\alpha^*).$$

The product of the associated orthogonal groups $O_\mathbf{d}=\prod_{i\in Q_0}O(W_i)$ is a maximal compact subgroup of ${\rm GL}_\mathbf{d}$ preserving the scalar product on $R_\mathbf{d}(Q)$. Denoting by $\bigoplus_{i\in Q_0}{\rm Symm}(W_i)$ the space of tuples of symmetric operators, we get the following description of the momentum map $\mu=\mu_Q:R_\mathbf{d}(Q)\rightarrow\bigoplus_{i\in Q_0}{\rm Symm}(W_i)$:
$$\mu((f_\alpha)_\alpha)=(\sum_{\alpha:~\rightarrow i}f_\alpha \circ f_\alpha^*-\sum_{\alpha:i\rightarrow ~}f_\alpha^*\circ f_\alpha)_i.$$

Then a direct application of \cite[Theorem 3.1]{OS} yields the following:

\begin{corollary}

Let $(B_i)_{i\in Q_0}$ be a tuple of symmetric operators in $\bigoplus_i{\rm Symm}(W_i)$, and let $a_{i,1}\leq\ldots\leq a_{i,d_i}$ be the eigenvalues of $B_i$, for all $i\in Q_0$. Then $(B_i)_i$ belongs to the image of $\mu_Q$ if and only if $$\sum_{i\in Q_0}\sum_k a_{i,k}=0\mbox{ and }a_\mathbf{K}=\sum_{i\in Q_0}\sum_{k\in K_i}a_{i,k}\geq 0\mbox{ for all }\mathbf{K}\in \hat{S}(\mathbf{d}),$$ where the sets $\hat{S}(\mathbf{d})$ are defined inductively as follows:\\[1ex]
a tuple of finite sets $\mathbf{K}$ belongs to $\hat{S}(\mathbf{d})$ if and only if
$${\rm ht}(\mathbf{K},\mathbf{L})\geq\sum_{\alpha:i\rightarrow j}|L_i|(d_j-|K_j|)\mbox{ for all }\mathbf{L}\in \hat{S}(|\mathbf{K}|).$$

\end{corollary}


\begin{thebibliography}{99}
\bibitem{A} M. F. Atiyah. Convexity and commuting Hamiltonians. Bull. London Math. Soc., 14(1):1–15, 1982.
\bibitem{BVW1} V. Baldoni, M. Vergne, M. Walter, Horn inequalities and quivers, arxiv:1804.00431
\bibitem{BVW2} V. Baldoni, M. Vergne, M. Walter, Horn conditions for Schubert positions of general quiver representations, arXiv:1901.07194
\bibitem{Br} M. Brion, Sur l’image de l’application moment, in: Seminaire d’algebre Paul Dubreil et Marie-Paule Malliavin, Proceedings, Paris 1986, Lecture Notes in Math. 1296 (Springer, Berlin, 1987), 177--192. 
\bibitem{CB} W. Crawley-Boevey, Geometry of the moment map for representations of quivers, Compositio Math., 126 (2001), 257--293.
\bibitem{CG} W. Crawley-Boevey, C. Geiss, Horn's problem and semi-stability for quiver representations, in: Representations of Algebras, Vol 1, Proceedings of the Ninth International Conference, Beijing, August 21-September 1, 2000, eds. D. Happel and Y. B. Zhang (Beijing Normal University Press, 2002), 40--48.
\bibitem{Ch} C. Chindris, Eigenvalues of Hermitian matrices and cones arising from quivers, International Mathematics Research Notices, 2006, Art. ID 59457, 27 pp.
\bibitem{GS} V. Guillemin and S. Sternberg. Convexity properties of the moment mapping, Invent.Math., 67(3):491–-513, 1982.
\bibitem{Kin} A. King, Moduli of representations of finite dimensional algebras,
Quarterly J. Math. Oxford 45 (1994) 515--530. 
\bibitem{Kir} F. Kirwan. Convexity properties of the moment mapping. III, Invent. Math.,77(3):547--552, 1984.
\bibitem{DW} H. Derksen and J. Weyman. Semi-Invariants of Quivers and Saturation for Littlewood-Richardson Coefficients, J. Amer. Math. Soc. 13 (2000), 467--479.
\bibitem{HLV} T. Hausel, E. Letellier and F. Rodriguez-Villegas, Positivity for Kac polynomials and DT-invariants of quivers, Ann. of Math.(2)177(2013), no. 3, 1147--1168.
\bibitem{Kac} V. G. Kac, Infinite Root Systems, Representations of Graphs and Invariant Theory, Inventiones math. 56 (1980), 57--92.
\bibitem{LBP} L. Le Bruyn, C. Procesi,  Semisimple representations of quivers, Trans. Amer. Math. Soc. 317 (1990), 585--598. 
\bibitem{OS} L. O'Shea, R. Sjamaar, Moment maps and Riemannian symmetric pairs, Math. Ann. 317 (2000), no. 3, 415-–457. 
\bibitem{Sch} A. Schofield,  General representations of quivers,  Proc. London Math. Socs. (3) 65 (1992), 1, 46--64.
\bibitem{Sj} R. Sjamaar, Convexity properties of the moment mapping re-examined, Adv. Math., 138(1):46-–91, 1998.



\end{thebibliography}
\end{document}